\documentstyle{article}\begin{document}\centerline{\bf Evaluation of an integral}\vskip .1in
\centerline{M.L. Glasser}
\centerline{Department of Theoretical Physics, University of Valladolid}

\centerline{Department of Physics,
Clarkson University}

\vskip .5in
\centerline{\bf  ABSTRACT}
\begin{quote}
The Moll-Arias de Reyna integral [1] 
$$\int_0^{\infty}\frac{dx}{(x^2+1)^{3/2}}\frac{1}{\sqrt{\varphi(x)+\sqrt{\varphi(x)}}}$$

$$\varphi(x)=1+\frac{4}{3}\left(\frac{x}{x^2+1}\right)^2$$
  is generalised and several values are given.

\end{quote}

\vskip .4in
\noindent
{\bf Keywwords}: Definite integral, Elliptic Integral

\vskip .3in
\noindent
{\bf 2010 MSC}: 33E05, 30E20

\newpage
\noindent
{\bf Introduction}\vskip .1in
 
We define{
$$f(a,b)=\int_0^{\infty}\frac{dx}{(x^2+1)^a}\frac{1}{\sqrt{\varphi(x)+\sqrt{\varphi(x)}}}\eqno(1)$$
where
$$\varphi(x)=1+4b^{-2}u^2,\quad u=\frac{x}{x^2+1}.\eqno(2)$$
The value $f(3/2, \sqrt{3})=\frac{\pi}{2\sqrt{6}}$ appeared as entry 3.248.5 in [2] and was shown to be incorrect by Moll et al.[3]. The exact value
$$f(3/2,\sqrt{3})=\frac{\sqrt{3}-1}{2}\Pi(\frac{\pi}{2},2-\sqrt{3},3^{-1/2})-6^{-1/2}F(\sin^{-1}\sqrt{2-\sqrt{3}},3^{-1/2})\eqno(3)$$
was recently provided in a mathematical {\it tour de force} by Arias de Reyna[1].
The aim of the present note is to provide further values of (1) and suggest that the incorrect value in [2] is not merely a misprint.
\vskip .1in\noindent
{\bf Calculation}\vskip .1in

By factoring $\sqrt{\varphi(x)}$  from the denominator of the integrand of (1), then multiplying the numerator and denominator by $\sqrt{\sqrt{\varphi(x)}-1}$ and changing the integration variable to $u$  (note that the range of integration wrt $x$ must first be divided into $[0,1]\cup[1,\infty]$) followed by $s=2u$ one obtains
$$f(a,b)=$$
$$2^{-a}b\int_0^1\frac{ds}{s\sqrt{1-s^2}}\left\{[1+\sqrt{1-s^2}]^{a-1}+[1-\sqrt{1-s^2}]^{a-1}\right\}\sqrt{1-\frac{b}{\sqrt{b^2+s^2}}}.\eqno(4)$$
Since both quadratic surds can be rationalised by the elliptic substitution
$s={\rm cn}(\kappa,x)$ for a suitable modulus, $f(a,b)$ should  be expressible in terms of elliptic integrals for integer and half integer values of $a$, even in the trigonometric case $\kappa=0$. This seems to eliminate the possibility of a simple misprint in [2]. For example, it is clear that
$$f(2,b)=\frac{1}{2}f(1,b)\eqno(5)$$
and with the substitution $t=b/\sqrt{b^2+s^2}$
$$f(1,b)=k\int_k^1\frac{dt}{(t+1)\sqrt{(1-t)(t^2-k^2)}},\quad k=\frac{b}{\sqrt{b^2+1}}.\eqno(6)$$
which is clearly a complete elliptic integral of the third kind and easily manipulated into standard form [4]
$$f(1,b)=\frac{k}{\sqrt{k+1}}\Pi(\frac{\pi}{2},\alpha^2,\kappa)$$
$$\alpha^2=\frac{\sqrt{b^2+1}+b}{2\sqrt{b^2+1}},\quad \kappa=\sqrt{b^2+1}-b.\eqno(7)$$

For $a=3,4,5,\cdots,$ $f(a,b)$, with $x=s^2$, can be easily seen to be a multiple of $f(1,b)$ plus an integral of the form
$$\int_0^1\frac{dx}{\sqrt{1-x}}P(x)\sqrt{1-\frac{b}{\sqrt{b^2+x}}}\eqno(8)$$
where $P$ is a polynomial having no constant term. Such an integral can always be manipulated into a sum of incomplete elliptic integrals  by the substitutions $x\rightarrow1-x^2$, $x\rightarrow\sqrt{b^2+1} \sin t$. For example,
$$f(3,b)
=\frac{1}{2}f(1,b)-\frac{k^2}{4}\int_{\frac{1}{\sqrt{1+k^2}}}^{1/k}\sqrt{\frac{x(x-1)}{1-k^2x^2}}dx.$$

For $a=3/2$ (4) yields
$$f(3/2,b)=\frac{b}{4}\int_0^{\pi/2}dt[{\rm csc}(t/2)+\sec(t/2)]\sqrt{1-\frac{b}{\sqrt{b^2+\sin^2t}}}.\eqno(10$$This can be further  simplified by the substitutions $\sin t=b\sin u$, $\sin u=x$, $\sqrt{1+x^2}=1/y$ to
$$f(3/2,\sqrt{3})=\frac{3}{\sqrt{8}}\sum_{\pm}\int_{\sqrt{3}/2}^1\frac{dy}{\sqrt{y(1+y)(4y^2-3)(y\pm\sqrt{4y^2-3})}}\eqno(11)$$
which may be reduced further to

$$f(2/3,\sqrt{3})=\frac{3^{1/4}}{2}\int_0^{1/\sqrt{3}}\frac{dx}{\sqrt{(x^2+1)(1-3x^2)}}\frac{\sqrt{X-2x}+\sqrt{X+2x}}{\sqrt{X(X+\frac{2}{\sqrt{3}})}}$$
$$ X=\sqrt{x^2+1}\eqno(12)$$
and which  may offer a more direct approach to (3).
\vskip .3in
\noindent
{\bf Acknowledgements}\vskip .1in
The author thanks Victor Moll for informing him of [1]. Partial financial support from Spanish Grants MTM2014-57129-C2-1-P (MINECO) and VA057U16 (Junta de Castilla y Le—n and FEDER) are gratefully acknowledged.

\newpage
\noindent
{\bf References}\vskip .2in
\noindent
[1] J. Arias de Reyna, ArXiv:1801.09640v1

\noindent
[2] I.S. Gradshteyn and I.M. Ryzhik, {\it Table of Integrals, Series and Products}. {Eds. A. Jeffrey, D. Zwillinger. Academic Press, New York, Sixth Edition, 2000].

\noindent
[3] T. Amdeberhan and V.H. Moll, {\it The integrals of Gradshteyn and Ryzhik. Part 14. An elementary evaluation of entry 2.411.5}, Sci. Ser. A.  Math. Sci. (N.S.) {\bf 19}(2010) 97-103.

\noindent
[4] P.F. Byrd and M.D. Friedman {\it Handbook of Elliptic Integrals fr Engineers and Scientists. Sec. Ed.} [Springer, Berlin (1971).]

\end{document}